\documentclass[12pt]{amsart}
\def\N{{\Bbb N}}
\def\Z{{\Bbb Z}}

\newtheorem{Theorem}{Theorem}[section]

\newtheorem{Lemma}[Theorem]{Lemma}

\theoremstyle{definition}

\newtheorem{Example}{Example}

\newtheorem*{Prob}{Problem}
\theoremstyle{remark}

\begin{document}
\sloppy
\title{A class of rigid Coxeter groups}
\author{Tetsuya Hosaka} 
\address{Department of Mathematics, Utsunomiya University, 
Utsunomiya, 321-8505, Japan}
\date{November 5, 2004}
\email{hosaka@cc.utsunomiya-u.ac.jp}
\thanks{Partly supported by the Grant-in-Aid for Scientific Research, 
The Ministry of Education, Culture, Sports, Science and Technology, Japan, 
(No.~15740029).}
\maketitle
\begin{abstract}
In this paper, 
we give a new class of rigid Coxeter groups.
Let $(W,S)$ be a Coxeter system. 
Suppose that 
(0) for each $s,t\in S$ such that $m(s,t)$ is even, $m(s,t)=2$, 
(1) for each $s\neq t\in S$ such that $m(s,t)$ is odd, 
$\{s,t\}$ is a maximal spherical subset of $S$, 
(2) there does not exist 
a three-points subset $\{s,t,u\}\subset S$ such that $m(s,t)$ and $m(t,u)$ are odd, 
and (3) for each $s\neq t\in S$ such that $m(s,t)$ is odd, 
the number of maximal spherical subsets of $S$ 
intersecting with $\{s,t\}$ is at most two, 
where $m(s,t)$ is the order of $st$ in the Coxeter group $W$.
Then we show that the Coxeter group $W$ is rigid.
This is an extension of a result of D.~Radcliffe.
\end{abstract}

\section{Introduction and preliminaries}

The purpose of this paper is to give a new class of rigid Coxeter groups.
A {\it Coxeter group} is a group $W$ having a presentation
$$\langle \,S \, | \, (st)^{m(s,t)}=1 \ \text{for}\ s,t \in S \,
\rangle,$$ 
where $S$ is a finite set and 
$m:S \times S \rightarrow \N \cup \{\infty\}$ is a function 
satisfying the following conditions:
(i) $m(s,t)=m(t,s)$ for any $s,t \in S$,
(ii) $m(s,s)=1$ for any $s \in S$, and
(iii) $m(s,t) \ge 2$ for any $s,t \in S$ such that $s\neq t$.
The pair $(W,S)$ is called a {\it Coxeter system}.
For a Coxeter group $W$, a generating set $S'$ of $W$ is called 
a {\it Coxeter generating set for $W$} if $(W,S')$ is a Coxeter system.
In a Coxeter system $(W,S)$, 
the conjugates of elements of $S$ are called {\it reflections}. 
Let $(W,S)$ be a Coxeter system.
For a subset $T \subset S$, 
$W_T$ is defined as the subgroup of $W$ generated by $T$, 
and called a {\it parabolic subgroup}.
If $T$ is the empty set, then $W_T$ is the trivial group.
A subset $T\subset S$ is called a {\it spherical subset of $S$}, 
if the parabolic subgroup $W_T$ is finite.

A {\it diagram} is an undirected graph $\Gamma$ 
without loops or multiple edges 
with a map $\text{Edges}(\Gamma)\rightarrow\{2,3,4,\ldots\}$ 
which assigns an integer greater than $1$ to each of its edges. 
Since such diagrams are used to define Coxeter systems, 
they are called {\it Coxeter diagrams}. 

Let $(W,S)$ and $(W',S')$ be Coxeter systems. 
Two Coxeter systems $(W,S)$ and $(W',S')$ are 
said to be {\it isomorphic}, 
if there exists a bijection 
$\psi:S\rightarrow S'$ such that 
$$m(s,t)=m'(\psi(s),\psi(t))$$ 
for every $s,t \in S$, where 
$m(s,t)$ and $m'(s',t')$ are the orders of $st$ in $W$ 
and $s't'$ in $W'$, respectively.

In general, a Coxeter group does not always determine 
its Coxeter system up to isomorphism.
Indeed some counter-examples are known.

\begin{Example}[{\cite[p.38 Exercise~8]{Bo}}, \cite{BMMN}]\label{Ex:B}
It is known that for an odd number $k\ge 3$, 
the Coxeter groups defined 
by the diagrams in Figure~\ref{fig1} are 
isomorphic and $D_{2k}$.
\begin{figure}[htbp]
\unitlength = 0.9mm
\begin{center}
\begin{picture}(80,28)(-40,-5)
\put(-25,0){\line(-1,0){20}}
\put(25,0){\line(1,0){20}}
\put(25,0){\line(2,3){10}}
\put(45,0){\line(-2,3){10}}
\put(-25,0){\circle*{1.3}}
\put(-45,0){\circle*{1.3}}
\put(25,0){\circle*{1.3}}
\put(45,0){\circle*{1.3}}
\put(35,15){\circle*{1.3}}
{\small
\put(-37,-4){$2k$}
\put(34,-4){$k$}
\put(27,8){$2$}
\put(41,8){$2$}
}
\end{picture}
\end{center}
\caption[Two distinct Coxeter diagrams for $D_{2k}$]{Two distinct Coxeter diagrams for $D_{2k}$}
\label{fig1}
\end{figure}
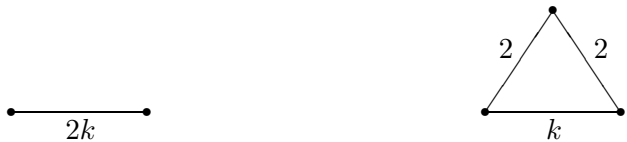
\end{Example}

\begin{Example}[\cite{BMMN}]\label{ex2}
It is known that
the Coxeter groups defined by 
the diagrams in Figure~\ref{fig2} are isomorphic by 
the {\it diagram twisting} (\cite[Definition~4.4]{BMMN}).
\begin{figure}[htbp]
\unitlength = 0.9mm
\begin{center}
\begin{picture}(140,30)(-70,-5)
\put(-5,0){\line(-1,0){60}}
\put(15,0){\line(1,0){40}}
\put(35,0){\line(0,1){20}}
\put(-5,0){\circle*{1.3}}
\put(-25,0){\circle*{1.3}}
\put(-45,0){\circle*{1.3}}
\put(-65,0){\circle*{1.3}}
\put(15,0){\circle*{1.3}}
\put(35,0){\circle*{1.3}}
\put(55,0){\circle*{1.3}}
\put(35,20){\circle*{1.3}}
{\small
\put(-16,-4){$2$}
\put(-36,-4){$3$}
\put(-56,-4){$2$}
\put(24,-4){$2$}
\put(44,-4){$2$}
\put(36,9){$3$}
}
\end{picture}
\end{center}
\caption[Coxeter diagrams for isomorphic Coxeter groups]{Coxeter diagrams for isomorphic Coxeter groups}
\label{fig2}
\end{figure}
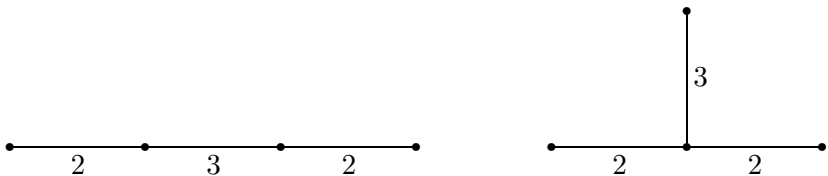
\end{Example}

Here there exists the following natural problem.

\begin{Prob}[\cite{BMMN}, \cite{CD}]
When does a Coxeter group determine its Coxeter system up to isomorphism?
\end{Prob}

A Coxeter group $W$ is said to be {\it rigid}, 
if for each Coxeter generating sets $S$ and $S'$ for $W$ 
the Coxeter systems $(W,S)$ and $(W,S')$ are isomorphic.

The following theorem was proved by D.~Radcliffe in \cite{R}.

\begin{Theorem}[\cite{R}, \cite{R2}]\label{Thm:R}
If $(W,S)$ is a right-angled Coxeter system 
(i.e.\ $m(s,t)=2$ for each $s\neq t\in S$), 
then the Coxeter group $W$ is rigid.
\end{Theorem}

In this paper, we prove the following theorem 
which is an extension of Theorem~\ref{Thm:R}.

\begin{Theorem}\label{Thm}
Let $(W,S)$ be a Coxeter system.
Suppose that 
\begin{enumerate}
\item[(0)] for each $s,t\in S$ such that $m(s,t)$ is even, 
$m(s,t)=2$, 
\item[(1)] for each $s\neq t\in S$ such that $m(s,t)$ is odd, 
$\{s,t\}$ is a maximal spherical subset of $S$, 
\item[(2)] there does not exist 
a three-points subset $\{s,t,u\}\subset S$ 
such that $m(s,t)$ and $m(t,u)$ are odd, and
\item[(3)] for each $s\neq t\in S$ such that $m(s,t)$ is odd, 
the number of maximal spherical subsets of $S$ 
intersecting with $\{s,t\}$ is at most two.
\end{enumerate}
Then the Coxeter group $W$ is rigid.
\end{Theorem}

Here we can not omit the condition (3) in this theorem by Example~\ref{ex2}.

\begin{Example}\label{ex3}
The Coxeter groups defined by 
the diagrams in Figure~\ref{fig3} are rigid by Theorem~\ref{Thm}.
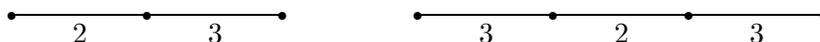
\begin{figure}[htbp]
\unitlength = 0.9mm
\begin{center}
\begin{picture}(140,10)(-70,-5)
\put(5,0){\line(1,0){60}}
\put(-15,0){\line(-1,0){40}}
\put(5,0){\circle*{1.3}}
\put(25,0){\circle*{1.3}}
\put(45,0){\circle*{1.3}}
\put(65,0){\circle*{1.3}}
\put(-15,0){\circle*{1.3}}
\put(-35,0){\circle*{1.3}}
\put(-55,0){\circle*{1.3}}
{\small
\put(14,-4){$3$}
\put(34,-4){$2$}
\put(54,-4){$3$}
\put(-26,-4){$3$}
\put(-46,-4){$2$}
}
\end{picture}
\end{center}
\caption[Coxeter diagrams for rigid Coxeter groups]{Coxeter diagrams for rigid Coxeter groups}
\label{fig3}
\end{figure}
\end{Example}

\section{Lemmas on Coxeter groups}

We first recall some basic results about Coxeter groups.

\begin{Lemma}[\cite{Bo}]\label{lem2-0}
Let $(W,S)$ be a Coxeter system and let $T\subset S$.
If there exists $w\in W$ such that $\ell(tw)<\ell(w)$ for every $t\in T$, 
then $T$ is a spherical subset of $S$.
Here $\ell(w)$ is the length of $w$ with respect to $S$.
\end{Lemma}

\begin{Lemma}[{\cite[p.12 Proposition~3]{Bo}}]\label{lem2-1}
Let $(W,S)$ be a Coxeter system and let $s,t\in S$.
Then $s$ is conjugate to $t$ if and only if 
there exists a sequence $s_1,\dots,s_n\in S$ such that 
$s_1=s$, $s_n=t$ and $m(s_i,s_{i+1})$ is odd for each $i\in\{1,\dots,n-1\}$.
\end{Lemma}

\begin{Lemma}[cf.\ \cite{BMMN}]\label{lem2-2}
Let $(W,S)$ and $(W',S')$ be Coxeter systems. 
Suppose that there exists an isomorphism $\phi:W\rightarrow W'$. 
Then for each maximal spherical subset $T\subset S$, 
there exists a unique maximal spherical subset $T'\subset S'$ 
such that $\phi(W_T)=w'W'_{T'}(w')^{-1}$ for some $w'\in W'$.
\end{Lemma}

Let $(W,S)$ be a Coxeter system and let $w\in W$.
A representation $w=s_1\cdots s_l$ ($s_i\in S$) 
is said to be {\it reduced}, if $\ell(w)=l$.

Using the above lemmas, we prove some lemmas needed later.

\begin{Lemma}\label{lem2-3-0}
Let $(W,S)$ be a Coxeter system which satisfies 
the conditions $(0)$ and $(1)$ in Theorem~\ref{Thm}, 
let $w,v,x\in W$ and 
let $w=s_1\cdots s_l$ and $v=t_1\cdots t_m$ ($s_i,t_j\in S$) 
be reduced representations.
Suppose that $W_{\{s_1,\dots,s_l\}}\cong \Z_2^l$, 
$W_{\{t_1,\dots,t_m\}}\cong \Z_2^m$ and $l\ge 2$.
If $w=xvx^{-1}$ then $w=v$.
\end{Lemma}

\begin{proof}
We prove this lemma by the induction on the length of $x$.
If $\ell(x)=0$ then $x=1$ 
and we obtain that if $w=xvx^{-1}$ then $w=v$.
Now we suppose that 
if $w=xvx^{-1}$ then $w=v$ 
for each $x\in W$ such that $\ell(x)<k$.
Let $x\in W$ such that $\ell(x)=k$.
Suppose that $w=xvx^{-1}$.
Let $x=a_1\cdots a_k$ ($a_i\in S$) be a reduced representation.
Then 
$$ s_1\cdots s_l=(a_1\cdots a_k)(t_1\cdots t_m)(a_k\cdots a_1). $$
If $(a_1\cdots a_k)(t_1\cdots t_m)$ is not reduced, 
then we may suppose that $a_k=t_1$, and 
$$ (a_1\cdots a_k)(t_1\cdots t_m)(a_k\cdots a_1)=
(a_1\cdots a_{k-1})(t_1\cdots t_m)(a_{k-1}\cdots a_1), $$
because $t_it_j=t_jt_i$ for each $i,j$.
Then by the inductive hypothesis, we obtain $w=v$.
By the same argument, 
if $(s_1\cdots s_l)(a_1\cdots a_k)$ is not reduced, 
then we obtain $w=v$ by the inductive hypothesis.
We suppose that 
$(a_1\cdots a_k)(t_1\cdots t_m)$ 
and $(s_1\cdots s_l)(a_1\cdots a_k)$ are reduced.
Here 
$$(s_1\cdots s_l)(a_1\cdots a_k)=(a_1\cdots a_k)(t_1\cdots t_m).$$
This implies that $l=m$.
Let $z=(s_1\cdots s_l)(a_1\cdots a_k)$.
Then $\ell(a_1z)<\ell(z)$ and 
$\ell(s_iz)<\ell(z)$ for each $i=1,\dots,l$.
By Lemma~\ref{lem2-0}, 
$\{s_1,\dots,s_l,a_1\}$ is a spherical subset of $S$.
Since $l\ge 2$, 
$W_{\{s_1,\dots,s_l,a_1\}}\cong \Z_2^{l+1}$ by (0) and (1).
Hence $a_1s_i=s_ia_1$ for each $i=1,\dots,l$ and 
$$ a_1z=(s_1\cdots s_l)(a_2\cdots a_k)=(a_2\cdots a_k)(t_1\cdots t_m),$$
i.e., $w=(a_2\cdots a_k)v(a_k\cdots a_2)$.
By the inductive hypothesis, we obtain $w=v$.
\end{proof}

\begin{Lemma}\label{lem2-3}
Let $(W,S)$ be a Coxeter system 
which satisfies the conditions $(0)$, $(1)$, $(2)$ and $(3)$ in Theorem~\ref{Thm}, 
let $x\in W$ and 
let $T$ and $U$ be subsets of $S$ 
such that $W_T\cong \Z_2^{|T|}$, $W_U\cong \Z_2^{|U|}$, 
$|T|\ge 2$ and $|U|\ge 2$.
Then $W_T\cap xW_Ux^{-1}\subset W_{T\cap U}$.
\end{Lemma}

\begin{proof}
Let $w\in W_T\cap xW_Ux^{-1}$.
Then $w=xvx^{-1}$ for some $v\in W_U$.
If $\ell(w)\ge 2$ then by Lemma~\ref{lem2-3-0}, 
$w=v\in W_T\cap W_U=W_{T\cap U}$.
We suppose that $\ell(w)=1$, i.e., $w\in T$. 
Let $t=w$.
Then $t=xvx^{-1}$ is a reflection and $v\in W_U$.
Hence $t=xvx^{-1}=ysy^{-1}$ for some $s\in U$ and $y\in W$, 
i.e., $t\sim s$.
Suppose that $t\neq s$.
Then $m(t,s)$ is odd by Lemma~\ref{lem2-1} and (2).
There exist maximal spherical subsets $\bar{T}$ and $\bar{U}$ of $S$ 
such that $T\subset\bar{T}$ and $U\subset\bar{U}$.
Here we note that $|\bar{T}|\ge|T|\ge 2$, $|\bar{U}|\ge|U|\ge 2$ and 
$\bar{T}\cap\bar{U}=\emptyset$ by (1).
Then $\{t,s\}$, $\bar{T}$ and $\bar{U}$ are different maximal spherical subsets 
intersecting with $\{t,s\}$.
This contradicts (3).
Thus $t=s$ and $w=t=s\in T\cap U\subset W_{T\cap U}$.
Hence $W_T\cap xW_Ux^{-1}\subset W_{T\cap U}$.
\end{proof}

\begin{Lemma}\label{lem2-4}
Let $(W,S)$ and $(W',S')$ be Coxeter systems. 
Suppose that there exists an isomorphism $\phi:W\rightarrow W'$ 
and that 
\begin{enumerate}
\item[$(2')$] there does not exist 
a three-points subset $\{s',t',u'\}\subset S'$ 
such that $m'(s',t')$ and $m'(t',u')$ are odd.
\end{enumerate}
Let $\{s,t\}$ and $T$ be maximal spherical subsets of $S$ 
such that $m(s,t)$ is odd.
By Lemma~\ref{lem2-2}, 
there exist maximal spherical subsets 
$\{s',t'\}$ and $T'$ of $S'$ 
such that $\phi(W_{\{s,t\}})\sim W'_{\{s',t'\}}$ and $\phi(W_T)\sim W'_{T'}$.
If $\{s,t\}\cap T\neq \emptyset$ then $\{s',t'\}\cap T'\neq \emptyset$.
\end{Lemma}

\begin{proof}
Suppose that $\{s,t\}\cap T\neq \emptyset$.
We may suppose that $t\in \{s,t\}\cap T$.
Then $\phi(t)\in\phi(W_{\{s,t\}})\sim W'_{\{s',t'\}}$.
Since $m(s,t)=m'(s',t')$ is odd and $o(\phi(t))=o(t)=2$, 
$\phi(t)$ is a reflection in $W'$ and $\phi(t)\sim s'\sim t'$.
Here $\phi(t)=x't'{x'}^{-1}$ for some $x'\in W'$ and 
$\phi(W_T)=w'W'_{T'}{w'}^{-1}$ for some $w'\in W'$.
Then 
$$ x't'{x'}^{-1}=\phi(t)\in \phi(W_T)=w'W'_{T'}{w'}^{-1}. $$
Hence $({w'}^{-1}x')t'({x'}^{-1}w')\in W'_{T'}$.
We note that $({w'}^{-1}x')t'({x'}^{-1}w')$ is a reflection in $W'_{T'}$.
Thus there exists $u'\in T'$ such that $t'\sim u'$.
Since $s'\sim t'\sim u'$, 
we have that $u'\in \{s',t'\}$
by the hypothesis $(2')$ and Lemma~\ref{lem2-1}.
Hence $\{s',t'\}\cap T'\neq \emptyset$.
\end{proof}

\section{Proof of the theorem}

Let $(W,S)$ be a Coxeter system.
Suppose that 
\begin{enumerate}
\item[(0)] for each $s,t\in S$ such that $m(s,t)$ is even, 
$m(s,t)=2$, 
\item[(1)] for each $s\neq t\in S$ such that $m(s,t)$ is odd, 
$\{s,t\}$ is a maximal spherical subset of $S$, 
\item[(2)] there does not exist 
a three-points subset $\{s,t,u\}\subset S$ 
such that $m(s,t)$ and $m(t,u)$ are odd, and
\item[(3)] for each $s\neq t\in S$ such that $m(s,t)$ is odd, 
the number of maximal spherical subsets of $S$ 
intersecting with $\{s,t\}$ is at most two.
\end{enumerate}
Let $(W',S')$ be a Coxeter system.
We suppose that there exists an isomorphism $\phi:W\rightarrow W'$.
The purpose is to show that the Coxeter systems 
$(W,S)$ and $(W',S')$ are isomorphic.

\begin{Lemma}\label{lem3-1}
The Coxeter system $(W',S')$ satisfies the following:
\begin{enumerate}
\item[$(0')$] for each $s',t'\in S'$ such that $m'(s',t')$ is even, 
$m'(s',t')=2$, 
\item[$(1')$] for each $s'\neq t'\in S'$ such that $m'(s',t')$ is odd, 
$\{s',t'\}$ is a maximal spherical subset of $S'$, 
\item[$(2')$] there does not exist 
a three-points subset $\{s',t',u'\}\subset S'$ 
such that $m'(s',t')$ and $m'(t',u')$ are odd, and
\item[$(3')$] for each $s'\neq t'\in S'$ such that $m'(s',t')$ is odd, 
the number of maximal spherical subsets of $S'$ 
intersecting with $\{s',t'\}$ is at most two.
\end{enumerate}
\end{Lemma}

\begin{proof}
$(0')$ and $(1')$ 
Let $s'\neq t'\in S'$.
There exists a maximal spherical subset $T'$ of $S'$ 
such that $\{s',t'\}\subset T'$.
By Lemma~\ref{lem2-2}, 
$\phi^{-1}(W'_{T'})\sim W_T$ for some maximal spherical subset $T$ of $S$.
By $(0)$ and $(1)$, either 
\begin{enumerate}
\item[(i)] $W_T\cong\Z_2^{|T|}$, or 
\item[(ii)] $|T|=2$ and if $T=\{s,t\}$ then $m(s,t)$ is odd.
\end{enumerate}
Hence $W_T$ is a rigid Coxeter group and 
$(W_T,T)$ and $(W'_{T'},T')$ are isomorphic.
Thus if $m'(s',t')$ is even then $m'(s',t')=2$, 
and if $m'(s',t')$ is odd then 
$\{s',t'\}$ is a maximal spherical subset of $S'$.
Hence $(0')$ and $(1')$ hold.

$(2')$ 
Suppose that there exists 
a three-points subset $\{s',t',u'\}\subset S'$ 
such that $m'(s',t')$ and $m'(t',u')$ are odd.
Then $\{s',t'\}$ and $\{t',u'\}$ are maximal spherical subsets of $S'$ by $(1')$.
By Lemma~\ref{lem2-2}, 
there exist $s,t,u,v\in S$ 
such that $\phi^{-1}(W'_{\{s',t'\}})\sim W_{\{s,t\}}$ 
and $\phi^{-1}(W'_{\{t',u'\}})\sim W_{\{u,v\}}$, 
where $\{s,t\}\neq\{u,v\}$ (i.e., $|\{s,t,u,v\}|\ge 3$).
We note that for each $x'\in W'_{\{s',t'\}}$, 
$x'$ is a reflection on $(W',S')$ if and only if $o(x')=2$, 
because $m'(s',t')$ is odd.
Hence for each $x'\in W'_{\{s',t'\}}$, 
$x'$ is a reflection on $(W',S')$ if and only if 
$\phi^{-1}(x')$ is a reflection on $(W,S)$.
Thus 
$$ s\sim t \sim \phi^{-1}(t')\sim\phi^{-1}(u')\sim u \sim v.$$
This contradicts (2) by Lemma~\ref{lem2-1}, since $|\{s,t,u,v\}|\ge 3$.
Thus $(2')$ holds.

$(3')$ 
Let $\{s,t\}$ and $T$ be maximal spherical subsets of $S$ 
such that $m(s,t)$ is odd.
By Lemma~\ref{lem2-2}, 
there exist maximal spherical subsets 
$\{s',t'\}$ and $T'$ of $S'$ 
such that $\phi(W_{\{s,t\}})\sim W'_{\{s',t'\}}$ and $\phi(W_T)\sim W'_{T'}$.
Then 
by Lemma~\ref{lem2-4}, $(2)$ and $(2')$, 
$\{s,t\}\cap T\neq \emptyset$ 
if and only if $\{s',t'\}\cap T'\neq \emptyset$.
Hence $(3')$ holds by $(3)$.
\end{proof}

\begin{Lemma}\label{lem3-2}
Let $T_1,\dots,T_k$ be maximal spherical subsets of $S$ and 
let $T'_1,\dots,T'_k$ be maximal spherical subsets of $S'$ 
such that $\phi(W_{T_i})\sim W'_{T'_i}$ and 
$W_{T_i}\cong \Z_2^{|T_i|}$ for each $i$.
Then $|T_1\cap\dots\cap T_k|=|T'_1\cap\dots\cap T'_k|$.
\end{Lemma}

\begin{proof}
If $|T_{i_0}|=1$ for some $i_0$, 
then $|T_1\cap\dots\cap T_k|=0=|T'_1\cap\dots\cap T'_k|$ 
because $T_i$ and $T'_i$ are maximal spherical subsets.
We suppose that $|T_i|\ge 2$ for each $i$.
Now $\phi(W_{T_i})=w'_iW'_{T'_i}{w'}_i^{-1}$ for some $w'_i\in W'$, and 
$$ \phi(W_{T_1\cap\dots\cap T_k})
=(w'_1W'_{T'_1}{w'}_1^{-1})\cap\dots\cap (w'_kW'_{T'_k}{w'}_k^{-1}) 
\subset w'_1W'_{T'_1\cap\dots\cap T'_k}{w'}_1^{-1}$$
by Lemmas~\ref{lem2-3} and \ref{lem3-1}.
Hence 
$\phi(W_{T_1\cap\dots\cap T_k})\subset 
w'_1W'_{T'_1\cap\dots\cap T'_k}{w'}_1^{-1}$.
By applying the same argument to $\phi^{-1}:W'\rightarrow W$, 
$$\phi^{-1}(W'_{T'_1\cap\dots\cap T'_k})\subset w_1W_{T_1\cap\dots\cap T_k}w_1^{-1}$$
for some $w_1\in W$.
Thus 
$$\phi(W_{T_1\cap\dots\cap T_k})\subset w'_1W'_{T'_1\cap\dots\cap T'_k}{w'}_1^{-1}
\subset w'_1\phi(w_1)\phi(W_{T_1\cap\dots\cap T_k})\phi(w_1^{-1}){w'}_1^{-1}.$$
Hence $\phi(W_{T_1\cap\dots\cap T_k})=w'_1W'_{T'_1\cap\dots\cap T'_k}{w'}_1^{-1}$.
Since $W_{T_1\cap\dots\cap T_k}\cong \Z_2^{|T_1\cap\dots\cap T_k|}$ 
and $W'_{T'_1\cap\dots\cap T'_k}\cong \Z_2^{|T'_1\cap\dots\cap T'_k|}$, 
we obtain that $|T_1\cap\dots\cap T_k|=|T'_1\cap\dots\cap T'_k|$.
\end{proof}

\begin{proof}[Proof of Theorem~\ref{Thm}]
Let $T_1,\dots,T_n$ be the maximal spherical subsets of $S$ and 
let $T'_1,\dots,T'_n$ be the maximal spherical subsets of $S'$ 
such that $\phi(W_{T_i})\sim W'_{T'_i}$ for each $i$.
By (0) and (1), either 
$W_{T_i}\cong \Z_2^{|T_i|}$, or 
$T_i=\{s_i,t_i\}$ and $m(s_i,t_i)$ is odd.
By Lemma~\ref{lem2-4}, Lemma~\ref{lem3-2}, $(2)$, $(3)$, $(2')$ and $(3')$, 
there exists a bijection $\psi:S\rightarrow S'$ such that 
$$\psi(\bigcap_{i\in I}T_i)=\bigcap_{i\in I}T'_i$$
for each $I\subset\{1,\dots,n\}$.
Then the bijection $\psi:S\rightarrow S'$ 
induces an isomorphism between $(W,S)$ and $(W',S')$.
\end{proof}

\section{Remark}

The following theorem has been proved in \cite{Ho}.

\begin{Theorem}[\cite{Ho}]\label{Thm:Ho}
Every Coxeter system $(W,S)$ which satisfies the conditions 
$(1)$, $(2)$ and $(3)$ in Theorem~\ref{Thm} is reflection rigid.
\end{Theorem}

Radcliffe has proved the following theorem in \cite{R2}.

\begin{Theorem}[\cite{R2}]\label{Thm:R2}
Let $(W,S)$ be a Coxeter system such that 
$m(s,t)\in \{2\}\cup 4\N$ for each $s\neq t\in S$.
Then the Coxeter group $W$ is rigid.
\end{Theorem}

Here the following problem arises as an extension of 
Theorems~\ref{Thm} and \ref{Thm:R2}.

\begin{Prob}
Let $(W,S)$ be a Coxeter system.
Suppose that 
\begin{enumerate}
\item[$(\bar{0})$] for each $s,t\in S$ such that $m(s,t)$ is even, 
$m(s,t)\in \{2\}\cup 4\N$, 
\item[$(1)$] for each $s\neq t\in S$ such that $m(s,t)$ is odd, 
$\{s,t\}$ is a maximal spherical subset of $S$, 
\item[$(2)$] there does not exist 
a three-points subset $\{s,t,u\}\subset S$ 
such that $m(s,t)$ and $m(t,u)$ are odd, and
\item[$(3)$] for each $s\neq t\in S$ such that $m(s,t)$ is odd, 
the number of maximal spherical subsets of $S$ 
intersecting with $\{s,t\}$ is at most two.
\end{enumerate}
Is it the case that the Coxeter group $W$ is rigid?
\end{Prob}

%

%
\end{document}